\documentclass[AMA,STIX1COL]{WileyNJD-v2}

\usepackage{moreverb}

\newcommand\BibTeX{{\rmfamily B\kern-.05em \textsc{i\kern-.025em b}\kern-.08em
T\kern-.1667em\lower.7ex\hbox{E}\kern-.125emX}}

\articletype{reasearch article}%

\received{<day> <Month>, <year>}
\revised{<day> <Month>, <year>}
\accepted{<day> <Month>, <year>}

\begin{document}

\title{Stochastic Maximum Principle for Control System with Time-varying delay}

\author[1]{Yuhang Li}

\author[2]{Yuecai Han}

\authormark{AUTHOR ONE \textsc{et al}}

\address{\orgdiv{School of Mathematics}, \orgname{Jilin University}, \orgaddress{\state{Changchun}, \country{China}}}

\corres{Yuecai Han, School of Mathematics, Jilin University,
Changchun 130012, China. \email{hanyc@jlu.edu.cn}}

\presentaddress{Jilin
University, 2699 Qianjin street, Changchun 130012, Jilin Province, China}

\abstract[Abstract]{In this paper, we study the stochastic optimal control problem for control system with time-varying delay. The corresponding stochastic differential equation is a kind of  stochastic differential delay equation. We prove the existence and uniqueness of the solution of this equation. We obtain the stochastic maximum principle of the control system with time-varying delay by introducing a kind of generalized anticipated backward stochastic differential equations. We prove the existence and uniqueness of the solution of this adjoint equation. As an application,  the linear quadratic moving average control problem is investigated  to illustrate the main result.}


\maketitle

\footnotetext{\textbf{Abbreviations:} ANA, anti-nuclear antibodies; APC, antigen-presenting cells; IRF, interferon regulatory factor}

\section{Introduction}\label{sec1}
Let $\tau(t)$ be a increasing differential random function such that $\tau(t)\le t\,  a.s.$ We focus on the following optimal control problem,
\begin{align}\label{1.1}
\left\{\begin{array}{ll}
dX_t =b\Big(t,X_t ,X_{\tau(t)},u_t,u_{\tau(t)}\Big)dt+\sigma\Big(t,X_t ,X_{\tau(t)},u_t,u_{\tau(t)} \Big)dW_t,\qquad 0\le t\le T,
\\X_0=x,
\end{array}\right.
\end{align}
to minimize the cost function
\begin{align}\label{cos}
    J(u)=E\left[\int_0^T f\left(t,X_t,X_{\tau(t)},u_t,u_{\tau(t)}\right)dt+g(X_T)\right].
\end{align}
We call the state equation (\ref{1.1}) as stochastic differential equation (SDE in short) with time-varying delay, which is a kind of stochastic delay differential equations (SDDEs in short). SDDEs and its control problem have been popularly investigated  and used in many areas \cite{oksendal2000maximum,arriojas2007delayed,chen2010maximum,mao2013delay,shuaiqi2020optimal}, etc. There are also many works about optimal control with time-varying delay in engineering and some other areas \cite{jajarmi2017efficient,sabermahani2020fibonacci,lin2022delay,guo2022dual,hu2022finite}. But as we know, there are not any works about maximum principle for control system with time-varying delay.

To study the  control system with delay, anticipated backward stochastic differential equations (ABSDEs in short) are introduced. ABSDE is first investigated by Peng and Yang \cite{peng2009anticipated} with the following form:
\begin{equation*}
\begin{cases}-d Y_t=f\left(t, Y_t, Z_t, Y_{t+\delta(t)}, Z_{t+\zeta(t)}\right) d t-Z_t d W_t, & t \in[0, T] , \\ Y_t=\xi_t, & t \in[T, T+K], \\ Z_t=\eta_t, & t \in[T, T+K].\end{cases}
\end{equation*}
Yang and Elliott \cite{yang2013some} study the properties of  generalized ABSDEs, they show the existence
and uniqueness of the solution and the duality between SDDEs. Hu et al. \cite{hu2021anticipated} extend the ABSDEs with quadratic growth, they study the solvability of these ABSDEs under different conditions. 
Thanks to the ABSDEs,  optimal control problem with delay is studied popularly. Chen and Wu\cite{chen2010maximum} obtain a maximum principle for control system with delay and apply it to a production and consumption choice problem.  Boccia and Vinter \cite{boccia2017maximum} provide different versions of the necessary conditions of optimality for optimal control problems
with delays, where  the conditions improve in a number of respects. Zhang et al.\cite{zhang2020stochastic} develop a stochastic maximum principle for partially-observed optimal control problems
for the state governed by stochastic differential equation with delay. Meng and Shi\cite{meng2021global} obtain the maximum principle where the control domain is non-convex and the diffusion term contains both control and its delayed
term.

Inspired by existing works, we study the maximum principle for optimal control problem with time-varying delay (\ref{1.1}), (\ref{cos}) in this paper. Here time-varying function $\tau(t)$  is independent with the Brownian motion $W(t)$. We study the uniqueness of the solution of the corresponding SDE with time-varying delay. Different from classical condition, to deal with the random delay, we use Gronwall inequality to $\sup_{0\le r\le t} E|\tilde{X}_r-X_r|^2$, which is a upper bound for   $E\left|\tilde{X}_{\tau(t)}-X_{\tau(t)}\right|^2$ under the independence assumption. Then a new stochastic maximum principle for control system (\ref{1.1}),(\ref{cos}) is established. We introduce the Hamiltonian function 
and the adjoint equation to obtain the optimal system.
To study the properties of the adjoint equation, we introduce the following ABSDE,
\begin{align*}
\left\{\begin{array}{ll}
-dy_t=h\left(t,y_t,z_t\right)dt+E^{\mathcal{F}_t}\left\{l\left(t,\theta(t),y_{\theta(t)},z_{\theta(t)}\right)d\theta(t)\right\}-z_tdW_t,\quad 0\le t\le T,
\\y_T=\xi ,
\\y_t=z_t=0,\quad t>T,
\end{array}\right.
\end{align*}
We prove the existence and uniqueness of it. Compared with classical type investigated by El Karoui and Peng\cite{el1997backward}, we construct a contraction mapping under a new $\beta$-norm 
\begin{align*}
\|(Y,Z)\|^2_{\beta}=\sup_{0\le s\le T}Ee^{\beta s}| Y_s|^2+ E\int_0^Te^{\beta s} Z_s^2ds.
\end{align*}
Furthermore, we get the  necessary condition that the optimal control process should satisfy.
Consider the linear quadratic
case,  we obtain the unique optimal control process for linear quadratic control system with time varying delay.

The rest of this paper is organized as follows. In section 2, we introduce the SDEs with time-varying delay and prove the existence and the uniqueness of the solution of this type of SDEs. 
In section 3, we obtain the stochastic maximum principle by introducing a kind of anticipated backward stochastic differential equations, and the existence and uniquenes of this kind of equations is proved.  
In the section 4, the linear quadratic case is investigated to illustrate the main results.

\section{Stochastic differential equation with time-varying delay}

\qquad Let $(\Omega,\mathcal{F},\mathbb{P})$ be a probability space,    $\mathcal{F}_0\subset \mathcal{F}$ be a sub $\sigma$-algebra. Let  $\mathbb{F}^W=(\mathcal{F}^W_t)_{0\le t\le T}$ be the filtration generated by a $m$-dimensional standard Brownian motion $\textbf{W}=(W_t)_{0\le t\le T}$ and  $\mathbb{F}^\tau=(\mathcal{F}^\tau_t)_{0\le t\le T}$ be the filtration generated by the process $\mathbb{\tau}=(\tau(t))_{0\le t\le T}$. Denote $\mathcal{F}_t=\mathcal{F}_0\lor\mathcal{F}^W_t\lor\mathcal{F}^\tau_t$ .  Consider the following SDE with time-varying delay,

\begin{align}\label{masde}
\left\{\begin{array}{ll}
dX_t =b(t,X_t ,X_{\tau(t)})dt+\sigma(t,X_t ,X_{\tau(t)})dW_t,\qquad 0\le t\le T,
\\X_0=x,
\end{array}\right.
\end{align}
where $E|X_0|^2< \infty$, $b$ and $\sigma $ be measurable functions on $[0,T]\times \mathbf{R}^d\times \mathbf{R}^d$ with values in $\mathbf{R}^d$ and $\mathbf{R}^{d\times m}$, respectively. 

Assume that 
\begin{equation}\label{2.3}
    |b(t,x,y)|^2\lor|\sigma (t,x,y)|^2 \le L(1+|x|^2+|y|^2),\quad x,y\in \mathbf{R}^n,t\in [0,T],
\end{equation}
and
\begin{align}\label{2.4}
|b(t,x_1,y_1)-b(t,x_2,y_2)|^2&\lor|\sigma(t,x_1,y_1)-\sigma(t,x_2,y_2)|^2
\le L(|x_1-x_2|^2+|y_1-y_2|^2),\notag\\ &x_1,x_2,y_1,y_2\in \mathbf{R}^n,t\in [0,T]
\end{align}
for some constant $L>0$ (where $|\sigma|^2=\sum |\sigma_{ij}|^2$).

\qquad Now we show the existence and uniqueness of the solution of equation (\ref{masde})
~\\

\textbf{Lemma 2.1}\qquad Assume that condition (\ref{2.3}) and (\ref{2.4}) holds, then there exist a unique solution to equation (\ref{masde}).
~\\

\textbf{\emph{Proof}:}
Uniqueness. Let $X_t$ and $\tilde{X}_t$ be two solutions of the equation (\ref{masde}),  we have

\begin{align*}
E|\tilde{X}_t-X_t|^2&=E\Big[\int_0^t b(s,\tilde{X}_s,\tilde{X}_{\tau(s)})-b(s,X_s,X_{\tau(s)})ds+\int_0^t \sigma(s,\tilde{X}_s,\tilde{X}_{\tau(s)})-\sigma(s,X_s,X_{\tau(s)})dW_s \Big]^2\notag\\
&\le2(T+1)LE\int_0^t |\tilde{X}_s-X_s|^2+|\tilde{X}_{\tau(s)}-X_{\tau(s)}|^2 ds\notag\\
&\le4(T+1)L\int_0^t\sup_{0\le r\le s}E|\tilde{X}_r-X_r|^2ds
.\end{align*}

For all $\varepsilon>0$, there exits $\xi_t\in[0,t]$, such that 
\begin{align*}
    E|\tilde{X}_{\xi_t}-X_{\xi_t}|^2\ge\sup_{0\le r\le t} E|\tilde{X}_r-X_r|^2-\varepsilon,
\end{align*}
so that
\begin{align*}
\sup_{0\le r\le t} E|\tilde{X}_r-X_r|^2\le&
    E|\tilde{X}_{\xi_t}-X_{\xi_t}|^2+\varepsilon\\
    \le&4(T+1)L\int_0^{\xi_t}\sup_{0\le r\le s}E|\tilde{X}_r-X_r|^2ds+\varepsilon\\
    \le&4(T+1)L\int_0^t\sup_{0\le r\le s}E|\tilde{X}_r-X_r|^2ds+\varepsilon
.\end{align*}
Through the Gronwall's inequality and the arbitrariness of $\varepsilon$, we get $\sup_{0\le t\le T} E|\tilde{X}_t-X_t|^2=0$. Thus, the solution $X_t$ is unique.
~\\

Existence. Let
\begin{align*}
\left\{\begin{array}{ll}
X_t^{(k+1)}&=x+\int_0^t b(s,X_s^{(k)},X_{\tau(s)}^{(k)})dt+\int_0^t\sigma(s,X_s^{(k)},X_{\tau(s)}^{(k)})dW_s ,
\\\quad X_t^{(0)}&=x,
\end{array}\right.
\end{align*}
and
\begin{align*}
    u_t^{(k)}=\sup_{0\le r\le t}E\Big|X_r^{(k+1)}-X_r^{(k)}\Big|^2
.\end{align*}
 Similar to the proof of classical case, we get

\begin{align*}
u_t^{(k)}\le\frac{A^{k+1}t^{k+1}}{(k+1)!}
\end{align*}
for some constants $A>0$. Let $\lambda$ be Lebesgue measure on $[0,T]$, $0\le n < m$ and $m,n\to \infty$. Then  we have
~\\

\begin{align*}
\left\|X_{t}^{(m)}-X_{t}^{(n)}\right\|_{L^2(\lambda\times P)}\le\sum_{k=n}^{m-1}\left(\frac{A^{k+2} T^{k+2}}{(k+2) !}\right)^{\frac{1}{2}} \rightarrow 0 
.\end{align*}
Therefore, $\{X_t^{(n)}\}_{n\ge 0}$ is a Cauchy sequence in ${L^2(\lambda\times P)}$. Define
\begin{align*}
    X_t:=\lim_{n\to\infty} X_t^{(n)}.
\end{align*}
Then $X_t$ is $\mathcal{F}_t$-measurable for all $t$. Since this holds for each $X_t^{(n)}$, thus $X_t$ is the solution of (\ref{masde}).

\section{The maximum principle}

\qquad Consider the following control problem. The state equation is
\begin{align}
\left\{\begin{array}{ll}
dX_t =b\Big(t,X_t ,X_{\tau(t)},u_t,u_{\tau(t)} \Big)dt+\sigma\Big(t,X_t ,X_{\tau(t)},u_t,u_{\tau(t)}  \Big)dW_t,\qquad 0\le t\le T,
\\X_0=x,
\end{array}\right.
\end{align}
with the cost function
\begin{align}
    J(u)=E\left[\int_0^T f\left(t,X_t ,X_{\tau(t)},u_t,u_{\tau(t)} \right)dt+g(X_T)\right],
\end{align}
where $b(t,x,y,u,v)$ and $\sigma(t,x,y,u,v) $ are measurable functions on $\mathbf{R}\times \mathbf{R}^d\times \mathbf{R}^d\times \mathbf{R}^k\times \mathbf{R}^k$ with values in $\mathbf{R}^d$ and $\mathbf{R}^{d\times m}$, respectively, and they satisfy Lipschitz's condition for $(x,y,u,v)$.
$f(t,x,y,u,v)$ and $g(x) $ be measurable functions on $\mathbf{R}\times \mathbf{R}^d\times \mathbf{R}^d\times \mathbf{R}^k\times \mathbf{R}^k$ and $\mathbf{R}^d$, respectively, with values in $\mathbf{R}$. 
We denote by $\mathbb{U}$ the set of control process
\textbf{u}$=(u_t)_{0\le t\le T}$ taking values in a given closed-convex set $\textbf{U}\subset \mathbb{R}^k$ and satisfying $E\int_0^T |u_t|^2 dt <\infty$.

To simplify the notation without losing the generality, we just consider the case $d=m=k=1$. Assume $u_t^{*}$ is the optimal control process, i.e.,
\begin{align*}
    J(u_t^*)=\inf_{u_t\in \mathbb{U}}J(u_t)
.\end{align*}
For $\forall\,0<\varepsilon<1$, let 
\begin{align*}
    u_t^\varepsilon =(1-\varepsilon)u_t^*+\varepsilon\alpha_t\triangleq 
u^*_t+\varepsilon v_t,
\end{align*}
where $\alpha_t$ is any other admissible control.

Define the Hamiltonian function $H$ by
\begin{align}
    H(t,x,y,u,v,p,q)=b(t,x,y,u,v)p+\sigma (t,x,y,u,v)q+f(t,x,y,u,v)
.\end{align}
Denote
\begin{align*}
l^*(t)=l\left(t,X_t^*,X^*_{\tau(t)},u_t^*,u^*_{\tau(t)}\right),
\end{align*}
for $l=b,\sigma,f,b_x,\sigma_x,f_x,b_y,\sigma_y,f_y,b_u,\sigma_u,f_u,b_v,\sigma_v,f_v$.
~\\

\textbf{Theorem 3.1} \qquad Assume that $(u_t^*)_{0\le t\le T}$ is the optimal control process, $(X_t^*)_{0\le t\le T}$  is the corresponding state process and $(p_t,q_t)$ are the processes satisfying
\begin{align}\label{absde}
\left\{\begin{array}{ll}
-dp_t=\left[b_x^*(t)p_t+\sigma_x^*(t)q_t+f_x^*(t)\right]dt+E^{\mathcal{F}_t}\left\{\left[b_y^*(\tau^{-1}(t))p_{\tau^{-1}(t)}+\sigma_y^*(\tau^{-1}(t))q_{\tau^{-1}(t)}+f_y^*(\tau^{-1}(t))\mathbf{1}_{\{\tau^{-1}(t)\le T\}}(t)\right]d\tau^{-1}(t)\right\}\\
\qquad\qquad
-q_tdW_t,\quad 0\le t\le T,
\\p_T=g_x(X_T^*) ,
\\p_t=q_t=0,\quad t>T.
\end{array}\right.
\end{align}
Then we have
\begin{align}\label{3.25}
    \left\{H^*_u(t)+E^{\mathcal{F}_t}\left[H_v^*\left(\tau^{-1}(t)\right)\mathbf{1}_{\{\tau^{-1}(t)\le T\}}(t)\right]\right\}\cdot(\alpha_t-u^*_t)\ge 0 ,
    \qquad \forall \alpha_t \in \mathbb{U},\quad \forall t\in[0,T],\quad a.s.
\end{align}
 for any control process $\alpha_t$, where 
\begin{align*}
    H^*(t)=H\Big(t,X_t^* ,X^*_{\tau(t)},u_t^* ,u^*_{\tau(t)},p_t,q_t\Big).
\end{align*}
~\\

\textbf{Remark 3.2} \quad Notice that $\tau^{-1}(t)$ is not $\mathcal{F}_t$-adapted, but $\{\tau^{-1}(t)\bigwedge T\}_{0\le t\le T}$ is $\mathcal{F}_T$-measurable. To investigate the adjoint equation (\ref{absde}), we consider a more general type of ABSDEs:
\begin{align}\label{aabsde}
\left\{\begin{array}{ll}
-dy_t=h\left(t,y_t,z_t\right)dt+E^{\mathcal{F}_t}\left\{l\left(t,\theta(t),y_{\theta(t)},z_{\theta(t)}\right)d\theta(t)\right\}-z_tdW_t,\quad 0\le t\le T,
\\y_T=\xi ,
\\y_t=z_t=0,\quad t>T,
\end{array}\right.
\end{align}
where $\{\theta(t)\}_{0\le t\le T}$ is a $\mathcal{F}_T$-measurable (not need to be adapted and differentiable) non-decreasing stochastic process such that $\theta(t)\ge t$ and $\theta(T)<\infty, a.s.$
For equation (\ref{absde}), we can take $\theta (t)=\tau^{-1}(t)\lor T$, which is bounded by T.

~\\

\textbf{Lemma 3.3} Assume that the following conditions hold:
\begin{align*}
|h(t,y,z)|\lor|l(t,\theta,y,z)|\le M_1(1+|y|+|z|)&,\\
|h(t,y_1,z_1)-h(t,y_2,z_2)|\lor|l(t,\theta,y_1,z_1)-l(t,\theta,y_2,z_2)|&\le M_1(|y_1-y_2|+|z_1-z_2|)
\end{align*}
and $\theta(T)<M_2, a.s.$ for some constants $M_1, M_2>0$ such that $16M_1M_2(1\lor M_2)<1$. 
Then the ABSDE (\ref{aabsde}) has the unique solution pair.

~\\

\textbf{\emph{Proof}:}
Denote $\mathbb{H}_T^2\left(\mathbb{R}^d\right)$ is the space of all predictable processes $\varphi: \Omega \times[0, T] \mapsto \mathbb{R}^d$ such that $\|\varphi\|^2=$ $\mathbb{E} \int_0^T\left|\varphi_t\right|^2 d t<+\infty.$

We define $\beta$-norm: $\|(Y,Z)\|^2_{\beta}=\sup_{0\le s\le T}Ee^{\beta s}| Y_s|^2+ E\int_0^Te^{\beta s} Z_s^2ds$ on $\mathbb{H}_T^2\left(\mathbb{R}^d\right)\times \mathbb{H}_T^2\left(\mathbb{R}^{d\times m}\right) $.
 For any $\mathcal{F}_t$-adapted  process pairs $(y_t^1,z_t^1), (y_t^2,z_t^2)$ with bounded $\beta$-norm and $y_t^1=y^2_t=z^1_t=z^2_t=0,\,t>T$, let
\begin{align}
\left\{\begin{array}{ll}
-dY^i_t=h\left(t,y^i_t,z^i_t\right)dt+E^{\mathcal{F}_t}\left\{l\left(t,\theta(t),y^i_{\theta(t)},z^i_{\theta(t)}\right)d\theta(t)\right\}-Z^i_tdW_t,\quad 0\le t\le T,
\\Y^i_T=\xi ,
\\Y^i_t=Z^i_t=0,\quad t>T
\end{array}\right.
\end{align}
for $i=1,2$.

Denote
\begin{align*}
\delta \phi_t=\phi^1_t-\phi^2_t
\end{align*}
for $\phi=Y,Z,y,z$. By it$\hat{\rm o}$'s formula,
\begin{align*}
d\left(e^{\beta t}\delta Y_t^2\right)=&\beta e^{\beta t}\delta Y_t^2dt+2e^{\beta t}\delta Y_td\delta Y_t+e^{\beta t}\left(d\delta Y_t\right)^2\notag\\
=&e^{\beta t}\left[\beta\delta Y_t^2-2\delta Y_t\delta h_t+\delta Z_t^2\right]dt-2e^{\beta t}E^{\mathcal{F}_t}[\delta Y_t\delta l_td\theta(t)]+e^{\beta t}\delta Y_t\delta Z_tdW_t,
\end{align*}
where $\delta h_t=h\left(t,y^1_t,z^1_t\right)-h\left(t,y^2_t,z^2_t\right)$ and $\delta l_t=l\left(t,\theta(t),y^1_{\theta(t)},z^1_{\theta(t)}\right)-l\left(t,\theta(t),y^2_{\theta(t)},z^2_{\theta(t)}\right)$.
Taking the  integral and expectation, we have
\begin{align*}
Ee^{\beta t}\delta Y_t^2+E\int_t^Te^{\beta s}\delta Z_s^2ds=&E\int_t^T e^{\beta s}\left[-\beta \delta Y_s^2+2\delta Y_s\delta h_s\right]ds+2E\int_t^T e^{\beta s}\delta Y_s\delta l_sd\theta(s)\notag\\
\le&E\int_t^Te^{\beta s}\left[-\beta\delta Y_s^2+c_1^{-1}\delta Y_s^2+c_1\delta h_s^2\right]ds+E\int_t^Te^{\beta s}\left[ c_2^{-1}\delta Y_s^2+c_2\delta l_s^2\right]d\theta(s)
\end{align*}
for all $c_1,c_2>0$. Notice that
\begin{align*}
E\int_t^Te^{\beta s}\delta Y_s^2d\theta (s)\le M_2\sup_{0\le s\le T}Ee^{\beta s}\delta Y_s^2,
\end{align*}
and 
\begin{align*}
E\int_t^T e^{\beta s}\delta l_s^2d\theta(s)\le M_1E\int_t^T e^{\beta \theta(s)}\left(\delta y_{\theta (s)}^2+\delta z_{\theta(s)}^2\right)d\theta(s)\le M_1M_2 \sup_{0\le s\le T}Ee^{\beta s}\delta y_s^2+M_1E\int_0^Te^{\beta s}\delta z_s^2ds,\quad \forall t\in[0,T].
\end{align*}
So we have
\begin{align*}
Ee^{\beta t}\delta Y_t^2+E\int_t^Te^{\beta s}\delta Z_s^2ds\le& E\int_t^Te^{\beta s}\left[(-\beta+c_1^{-1})\delta Y_s^2+c_1M_1(\delta y_s^2+\delta z_s^2)\right]ds+c_2^{-1}M_2\sup_{0\le s\le T}Ee^{\beta s}\delta Y_s^2\\
&+c_2M_1M_2 \sup_{0\le s\le T}Ee^{\beta s}\delta y_s^2+c_2M_1E\int_0^Te^{\beta s}\delta z_s^2ds,\quad \forall \in [0,T].
\end{align*}
Let $\beta>c_1^{-1}$ and $c_1\to 0$. We get
\begin{align*}
(1-c_2^{-1}M_2)\sup_{0\le s\le T}Ee^{\beta s}\delta Y_s^2\le c_2M_1M_2 \sup_{0\le s\le T}Ee^{\beta s}\delta y_s^2+c_2M_1E\int_0^Te^{\beta s}\delta z_s^2ds,
\end{align*}
and
\begin{align*}
E\int_0^Te^{\beta s}\delta Z_s^2ds\le c_2^{-1}M_2\sup_{0\le s\le T}Ee^{\beta s}\delta Y_s^2+
c_2M_1M_2 \sup_{0\le s\le T}Ee^{\beta s}\delta y_s^2+c_2M_1E\int_0^Te^{\beta s}\delta z_s^2ds.
\end{align*}
So we have
\begin{align*}
(1-2c_2^{-1}M_2)\sup_{0\le s\le T}Ee^{\beta s}\delta Y_s^2+E\int_0^Te^{\beta s}\delta Z_s^2ds\le 2c_2M_1M_2 \sup_{0\le s\le T}Ee^{\beta s}\delta y_s^2+2c_2M_1E\int_0^Te^{\beta s}\delta z_s^2ds,
\end{align*}
which shows
\begin{align*}
(1-2c_2^{-1}M_2)\|(\delta Y,\delta Z)\|^2_{\beta}\le2c_2M_1(1\lor M_2)\|(\delta y,\delta z)\|^2_{\beta}.
\end{align*}

Under the assumption $16M_1M_2(1\lor M_2)<1$ and taking $c_2=\left(\frac{M_2}{M_1(1\lor M_2)}\right)^{\frac{1}{2}}$, we get the contraction mapping $T:(y,z)\to(Y,Z)$ from $\mathbb{H}_{T,\beta}^2\left(\mathbb{R}^d\right)\times\mathbb{H}_{T,\beta}^2\left(\mathbb{R}^{d\times m}\right)$ 
onto itself, so that there exists a fixed point, which is the unique continuous solution of
the  ABSDE (\ref{aabsde}).

This completes the proof of Lemma 3.3.

~\\

\textbf{Lemma 3.4} \qquad Let $(u_t^*)_{0\le t\le T}$ be the optimal control process and $(X_t^*)_{0\le t\le T}$ be the corresponding state process, and $(p_t,q_t)$ is the adjoint process satisfying(\ref{absde}).
Then the $G\hat{a}teaux$ derivative of $J$ at $u^*_t$ in the direction $v_t$ is
\begin{align}
    \frac{d}{d\varepsilon}J(u_t^*+\varepsilon v_t)\Big|_{\varepsilon=0}=E\int_0^T \Big[H^*_u(t)+H^*_v(\tau^{-1}(t))\frac{\partial \tau^{-1}(t)}{\partial t}\mathbf{1}_{\{\tau^{-1}(t)<T\}}(t)\Big]\cdot v_t dt
.\end{align}

\textbf{\emph{Proof}:}\quad Let $X_t^*$ and $X_t^\varepsilon$ be the state process corresponding to $u_t^*$ and $u_t^\varepsilon$, respectively. Define $V_t$ by
\begin{align}
\left\{\begin{array}{ll}
dV_t=&\Big[b_x^*(t) V_t+b_y^*(t)V_{\tau(t)}+b_u^*(t) v_t+b_v^*(t)v_{\tau(t)}\Big]dt\\
&+\Big[\sigma_x^*(t) V_t+\sigma_y^*(t)V_{\tau(t)}+\sigma_u^*(t) v_t+\sigma_v^*(t)v_{\tau(t)}\Big]d W_t,
\\V_0=0.
\end{array}\right.
\end{align}
 Then it's easy to get
\begin{align*}
\sup_{0\le t\le T}\lim_{\varepsilon\to 0} E\Big[\frac{X_t^\varepsilon-X_t^*}{\varepsilon}-V_t\Big]^2=0.
\end{align*}
~\\
So we have                              
\begin{align}\label{3.7}
\frac{J(u_t^\varepsilon)-J(u_t^*)}{\varepsilon}\to
E\left[\int_0^T\left(f_x^*(t)V_t+f_y^*(t)V_{\tau(t)}+f_u^*(t)v_t+f_v^*(t)v_{\tau(t)}\right)dt+g_x(X_T^*)V_T\right]
,\end{align}
as $\varepsilon\to 0$.

~\\

 By It$\hat{\rm o}$'s formula, we have that
\begin{align}
d(p_tV_t)&=p_tdV_t+V_tdp_t+dp_tdV_t\notag \\
    &=p_t\Big[b_x^*(t) V_t+b_y^*(t)V_{\tau(t)}+b_u^*(t) v_t+b_v^*(t)v_{\tau(t)}\Big]dt\notag -V_t\left[b_x^*(t)p_t+\sigma_x^*(t)q_t+f_x^*(t)\right]dt\\
    &\quad -V_tE^{\mathcal{F}_t}\left\{\left[b_y^*(\tau^{-1}(t))p_{\tau^{-1}(t)}+\sigma_y^*(\tau^{-1}(t))q_{\tau^{-1}(t)}+f_y^*(\tau^{-1}(t))\mathbf{1}_{\{\tau^{-1}(t)\le T\}}(t)\right]d\tau^{-1}(t)\right\}\notag \\
    &\quad+q_t\Big[\sigma_x^*(t) V_t+\sigma_y^*(t)V_{\tau(t)}+\sigma_u^*(t) v_t+\sigma_v^*(t)v_{\tau(t)}\Big]dt+M_tdW_t \notag \\
    &=p_t\Big[b_y^*(t)V_{\tau(t)}+b_u^*(t) v_t+b_v^*(t)v_{\tau(t)}\Big]dt -f_x^*(t)V_tdt\notag\\
    &\quad -V_tE^{\mathcal{F}_t}\left\{\left[b_y^*(\tau^{-1}(t))p_{\tau^{-1}(t)}+\sigma_y^*(\tau^{-1}(t))q_{\tau^{-1}(t)}+f_y^*(\tau^{-1}(t))\mathbf{1}_{\{\tau^{-1}(t)\le T\}}(t)\right]d\tau^{-1}(t)\right\}\notag \\
    &\quad+q_t\Big[\sigma_y^*(t)V_{\tau(t)}+\sigma_u^*(t) v_t+\sigma_v^*(t)v_{\tau(t)}\Big]dt+M_tdW_t
,\end{align}
where $(M_t)_{0\le t\le T}$ is $\mathcal{F}_t$-adapted. Notice that
\begin{align}
\int_0^Tb^*_y(\tau^{-1}(t))p_{\tau^{-1}(t)}V_td\tau^{-1}(t)=\int_0^{\tau^{-1}(T)}b_y^*(t)p_tV_{\tau(t)}dt=\int_0^Tb_y^*(t)p_tV_{\tau(t)}dt
\end{align}
and
\begin{align}
\int_0^T\sigma^*_y(\tau^{-1}(t))q_{\tau^{-1}(t)}V_td\tau^{-1}(t)=\int_0^T\sigma_y^*(t)q_tV_{\tau(t)}dt.
\end{align}
Then
\begin{align}\label{3.14}
    Eg_x(X_T^*)V_T
    &=Ep_TV_T=E\int_0^T d(p_tV_t)\notag\\
    &=E\int_0^T \left[p_tb_u^*(t)v_t+p_tb_v^*(t)v_{\tau(t)}+q_t\sigma_u^*(t)v_t+q_t\sigma_v^*(t)v_{\tau(t)}\right]dt\notag\\
    &\quad -E\int_0^Tf_x^*(t)V_tdt-E\int_0^Tf_y^*(\tau^{-1}(t))V_t\mathbf{1}_{\{\tau^{-1}(t)\le T\}}d\tau^{-1}(t)
.\end{align}
So we have
\begin{align}
    \frac{d}{d\varepsilon}J(u_t^*+\varepsilon v_t)\Big|_{\varepsilon=0}&=E\int_0^T\Big(f_x^*V_t+f_y^*V_{\tau(t)}+f_u^*v_t+f_v^*(t)v_{\tau(t)}\Big)dt\notag\\
    &+E\int_0^T \left[p_tb_u^*(t)v_t+p_tb_v^*(t)v_{\tau(t)}+q_t\sigma_u^*(t)v_t+q_t\sigma_v^*(t)v_{\tau(t)}\right]dt\notag\\
    &\quad -E\int_0^Tf_x^*(t)V_tdt-E\int_0^Tf_y^*(\tau^{-1}(t))V_t\mathbf{1}_{\{\tau^{-1}(t)\le T\}}d\tau^{-1}(t)\\
    &=E\Big[\int_0^Tb_u^*(t)p_t+\sigma_t^*q_t+f_u^*(t)\Big]v_tdt+E\int_0^T\Big[b_v^*(t)p_t+\sigma_v^*(t)q_t+f_v^*(t)\Big]v_{\tau(t)}dt\notag\\
    &=E\int_0^T\Big[ H_u^*(t) v_t+H_v^*(t)v_{\tau(t)}\Big] dt\notag\\
    &=E\int_0^T \Big[H^*_u(t)+H^*_v(\tau^{-1}(t))\frac{\partial \tau^{-1}(t)}{\partial t}\mathbf{1}_{\{\tau^{-1}(t)<T\}}(t)\Big]\cdot v_t dt
\end{align}
 This completes the proof of Lemma 3.4.
~\\

Since $(u_t^*)_{0\le t\le T}$ is optimal control process, we have the inequality
\begin{align*}
    \frac{d}{d\varepsilon}J\Big(u_t^*+\varepsilon v_t\Big)\Big|_{\varepsilon=0}\ge 0
.\end{align*}
By Lemma 3.4, we get
\begin{align*}
   E\int_0^T \Big[H^*_u(t)+H^*_v(\tau^{-1}(t))\frac{\partial \tau^{-1}(t)}{\partial t}\mathbf{1}_{\{\tau^{-1}(t)<T\}}(t)\Big]\cdot v_t dt\ge 0
.\end{align*}
So
\begin{align*}
    E\Big[\mathbf{1}_A\big[H^*_u(t)+H^*_v(\tau^{-1}(t))\frac{\partial \tau^{-1}(t)}{\partial t}\mathbf{1}_{\{\tau^{-1}(t)<T\}}(t)\big]\Big]\cdot v_t\ge 0,\quad \forall t\in[0,T],\quad \forall A\subset \mathcal{F}_t.
\end{align*}
To ensure adaptability, we can rewrite the above equation as
\begin{align*}
    E\Big[\mathbf{1}_A\big[H^*_u(t)+E^{\mathcal{F}_t}[H^*_v(\tau^{-1}(t))\frac{\partial \tau^{-1}(t)}{\partial t}\mathbf{1}_{\{\tau^{-1}(t)<T\}}(t)]\big]\Big]\cdot v_t\ge 0,\quad \forall t\in[0,T],\quad \forall A\subset \mathcal{F}_t,
\end{align*}
and obtain that
\begin{align*}
    \Big[H^*_u(t)+E^{\mathcal{F}_t}[H^*_v(\tau^{-1}(t))\frac{\partial \tau^{-1}(t)}{\partial t}\mathbf{1}_{\{\tau^{-1}(t)<T\}}(t)]\Big]\cdot v_t \ge 0,\qquad \forall t\in[0,T].
\end{align*}
This completes the proof of Theorem 3.1.

~\\

\textbf{Remark 3.5}  \qquad If the optimal control process $(u_t^*)_{0\le t\le T}$ takes values in the interior of the $\mathbb{U}$ , then we can replace (\ref{3.25})  with the following condition
\begin{align*}
    H_u^*(t)+E^{\mathcal{F}_t} 
 \Big[H^*_v(\tau^{-1}(t))\frac{\partial \tau^{-1}(t)}{\partial t}\mathbf{1}_{\{\tau^{-1}(t)<T\}}(t)\Big]=0
.\end{align*}
~\\

Thus, we give the optimal system
 \begin{align*}
    \left\{\begin{array}{ll}
dX_t^*=H_p^*(t)dt+H_q^*(t)dW_t,\\\\
  -dp_t=H_x^*(t)dt+E^{\mathcal{F_t}}\left[H_y^*\left(\tau^{-1}(t)\right)\mathbf{1}_{\{\tau^{-1}(t)\le T\}}d\tau^{-1}(t)\right] -q_tdW_t,\\\\
  X_0^*=x,\quad  p_T=g_x(X_T^*),\\\\
H_u^*(t)+E^{\mathcal{F}_t} 
 \Big[H^*_v(\tau^{-1}(t))\frac{\partial \tau^{-1}(t)}{\partial t}\mathbf{1}_{\{\tau^{-1}(t)<T\}}(t)\Big]=0,
\end{array}\right.
\end{align*}

where
\begin{align*}
    H^*(t)&=H\Big(t,X_t^* ,X^*_{\tau(t)},u_t^* ,u^*_{\tau(t)},p_t,q_t\Big),\\
    H(t,x,y,u,v,p,q)&=b(t,x,y,u,v)p+\sigma(t,x,y,u,v)q+f(t,x,y,u,v).
\end{align*}

\section{Linear quadratic case}

\quad In this section, we consider a linear quadratic (LQ in short) case, which can describe a moving average linear quadratic regulator problem. The state process is defined as follows
\begin{align}\label{4.1}
    dX_t=\Big(A_tX_t+B_tX_{at}+C_tu_t+P_tu_{at}\Big)dt+\Big(D_tX_t+F_tX_{at}+H_tu_t+N_tu_{at}\Big)dW_t ,
\end{align}
with the cost function
\begin{align}\label{4.2}
    J(u)=\frac{1}{2}E\Big[\int_0^T (Q_tX_t^2+S_tX_{at}^2+R_tu_t^2)dt+GX_T^2\Big].
\end{align}
Here $a\in(0,1), G>0$ and $Q_t,S_t,R_t$ are positive definite functions.
~\\

Using the conclusions of Section 3, we can get the adjoint equation
\begin{align}\label{4.3}
    \left\{\begin{array}{ll}
-dp_t=\left[A_tp_t+D_tq_t+Q_t+a^{-1}\left[B_{a^{-1}t}E^{\mathcal{F}_t}[p_{a^{-1}t}]+F_{a^{-1}t}E^{\mathcal{F}_t}[q_{a^{-1}t}]+S_{a^{-1}t}\mathbf{1}_{[0,aT]}(t)\right]\right]dt\\
\\
 \qquad \qquad -q_tdW_t,\quad t\in[0,T]
 \\
\\p_T=GX_T^* ,
\\
\\p_t=q_t=0, \quad t>T,
\end{array}\right.
\end{align}
 and the optimal control process $u_t^*$ should satisfy
 \begin{align*}
     C_tp_t+H_tq_t+R_tu^*_t+a^{-1}P_{a^{-1}t}E^{\mathcal{F}_t}[p_{a^{-1}t}]+a^{-1}N_{a^{-1}t}E^{\mathcal{F}_t}[q_{a^{-1}t}]=0,
 \end{align*}
i.e.,
\begin{align}\label{4.4}
    u_t^*=-R_t^{-1}\left(C_tp_t+H_tq_t+a^{-1}P_{a^{-1}t}E^{\mathcal{F}_t}[p_{a^{-1}t}]+a^{-1}N_{a^{-1}t}E^{\mathcal{F}_t}[q_{a^{-1}t}]\right).
\end{align}

\textbf{Theorem 4.1} The function $ u_t^*=-R_t^{-1}\left(C_tp_t+H_tq_t+a^{-1}P_{a^{-1}t}E^{\mathcal{F}_t}[p_{a^{-1}t}]+a^{-1}N_{a^{-1}t}E^{\mathcal{F}_t}[q_{a^{-1}t}]\right),\quad t\in[0,T]$ is the unique optimal control for moving average LQ problem (\ref{4.1}), (\ref{4.2}), where $(p_t,q_t)$ is defined by equation (\ref{4.3}). 
~\\

\textbf{\emph{Proof}:} \quad We now prove $u_t^*$ is the optimal control. For any $\tilde{u}_t\subset \mathbb{U}$, let $(\tilde{X}_t,\tilde{X}_{at},\tilde{u}_{at})$ and $(X_t^*,X_{at}^*,u_{at}^*)$ are processes corresponding to $\tilde{u}_t$ and $u_t^*$, respectively. We have that 
\begin{align*}
    d(\tilde{X}_t-X_t^*)=&[A_t(\tilde{X}_t-X^*_t)+B_t(\tilde{X}_{at}-X^*_{at})+C_t(\tilde{u}_t-u^*_t)+P_t(\tilde{u}_{at}-u^*_{at})]dt\\
    &+[D_t(\tilde{X}_t-X^*_t)+F_t(\tilde{X}_{at}-X^*_{at})+H_t(\tilde{u}_t-u^*_t)+N_t(\tilde{u}_{at}-u^*_{at})]dW_t.
\end{align*}
Consider
\begin{align}\label{4.5}
    dp_t(\tilde{X}_t-X_t^*)=&p_td(\tilde{X}_t-X^*_t)+(\tilde{X}_t-X^*_t)dp_t+dp_td(\tilde{X}_t-X_t^*)\notag\\
    =&p_t[A_t(\tilde{X}_t-X^*_t)+B_t(\tilde{X}_{at}-X^*_{at})+C_t(\tilde{u}_t-u^*_t)+P_t(\tilde{u}_{at}-u^*_{at})]dt\notag\\
    &-(\tilde{X}_t-X_t^*)\left[A_tp_t+D_tq_t+Q_t+a^{-1}\left[B_{a^{-1}t}E^{\mathcal{F}_t}[p_{a^{-1}t}]+F_{a^{-1}t}E^{\mathcal{F}_t}[q_{a^{-1}t}]+S_{a^{-1}t}\mathbf{1}_{[0,aT]}(t)\right]\right]dt\notag\\ 
    &+q_t\left[D_t(\tilde{X}_t-X^*_t)+F_t(\tilde{X}_{at}-X^*_{at})+H_t(\tilde{u}_t-u^*_t)+N_t(\tilde{u}_{at}-u^*_{at})\right]dt+M_tdW_t\notag\\
    =&\left[B_tp_t(\tilde{X}_{at}-X^*_{at})-a^{-1}B_{a^{-1}t}E^{\mathcal{F}_t}[p_{a^{-1}t}]\right](\tilde{X}_t-X_t^*)dt\notag\\
    &+\left[F_tq_t(\tilde{X}_{at}-X^*_{at})-a^{-1}F_{a^{-1}t}E^{\mathcal{F}_t}[q_{a^{-1}t}]\right](\tilde{X}_t-X_t^*)dt\notag\\
    &+\left[C_tp_t+H_tq_t\right](\tilde{u}_t-u^*_t)dt+\left[P_tp_t+N_tq_t\right](\tilde{u}_{at}-u^*_{at})dt\notag\\
    &-\left[Q_t+a^{-1}S_{a^{-1}t}\mathbf{1}_{[0,aT]}(t)\right](\tilde{X}_t-X_t^*)dt+M_tdW_t,
\end{align}
where $(M_t)_{0\le t\le T}$ is a $\mathcal{F}_t$-adapted process.
Taking integral from $0$ to $T$ and taking the expectation and by (\ref{4.4}), we have
\begin{align}
EGX_T^*(\tilde{X}_T-X_T^*)=&Ep_T(\tilde{X}_T-X_T^*)\notag=E\int_0^Tdp_t(\tilde{X}_t-X_t^*)\notag\\
    =&-E\int_0^T\Big[R_tu_t^*(\tilde{u}_t-u^*_t)+Q_tX_t^*(\tilde{X}_t-X^*_t)+S_tX^*_{at}(\tilde{X}_{at}-X_{at}^*)\Big]dt.
\end{align}
By using the fact $a^2-b^2\ge 2b(a-b)$, we have that
\begin{align*}
J(\tilde{u}_t)-J(u_t^*)\ge&\frac{1}{2}E\int_0^T \Big[Q_t(\tilde{X}_t^2-X_t^{*2})+S_t(\tilde{Y}_t^2-Y_t^{*2})+R_t(\tilde{u}_t^2-u_t^{*2})\Big]dt \\
    &+\frac{1}{2}EGX^*_T(\tilde{X}_T-X_T^*) \\
    =& \frac{1}{2}E\int_0^T \Big[Q_t(\tilde{X}_t^2-X_t^{*2})-2Q_tX_t^*(\tilde{X}_t-X^*_t)+S_t(\tilde{X}_{at}^2-X_{at}^{*2})\\
     &\qquad\qquad-2S_tX^*_{at}(\tilde{X}_{at}-X_{at}^*)
   +R_t(\tilde{u}_t^2-u_t^{*2})-2R_tu_t^*(\tilde{u}_t-u^*_t)\Big]dt \\
    \ge& 0.
\end{align*}
This shows that $u_t^*$ is an optimal control. 

Then we prove $u^*_t$ is unique. Assume that both $u_t^{*,1}$ and $u_t^{*,2}$ are optimal controls, $X_t^1$ and $X_t^2$ are corresponding state processes, respectively. It is easy to get $\frac{X_t^1+X_t^2}{2}$ is the corresponding state process to $\frac{u_t^{*,1}+u_t^{*,2}}{2}$. We assume there exist constants $\delta>0, \alpha\ge 0$, such that $R_t\ge \delta$ and 
\begin{align*}
    J(u_t^{*,1})=J(u_t^{*,2})=\alpha.
\end{align*}
Using the fact $a^2+b^2=2[(\frac{a+b}{2})^2+(\frac{a-b}{2})^2]$, we have that

\begin{align*}
    2\alpha=&J(u_t^{*,1})+J(u_t^{*,2})\\
           =&\frac{1}{2}E\int_0^T \Big[Q_t(X_t^1X_t^1+X_t^2X_t^2)+S_t(X_{at}^1X_{at}^1+X_{at}^2X_{at}^2)+R_t(u_t^{*,1}u_t^{*,1}+u_t^{*,2}u_t^{*,2})\Big]dt\\
           &+\frac{1}{2}EG(X_T^1X_T^1+X_T^2X_T^2)\\
           \ge&E\int_0^T \Big[Q_t\Big(\frac{X_t^1+X_t^2}{2}\Big)^2+S_t\Big(\frac{X_{at}^1+X_{at}^2}{2}\Big)^2+R_t\Big(\frac{u_t^{*,1}+u_t^{*,2}}{2}\Big)^2\Big]dt\\
           &+EG\Big(\frac{X_T^1+X_T^2}{2}\Big)^2+E\int_0^TR_t\Big(\frac{u_t^{*,1}-u_t^{*,2}}{2}\Big)^2dt\\
           =&2J\Big(\frac{u_t^{*,1}+u_t^{*,2}}{2}\Big)+E\int_0^TR_t\Big(\frac{u_t^{*,1}-u_t^{*,2}}{2}\Big)^2dt\\
           \ge&2\alpha+\frac{\delta}{4}E\int_0^T|u_t^{*,1}-u_t^{*,2}|^2dt.
\end{align*}
Thus, we have 
\begin{align*}
    E\int_0^T|u_t^{*,1}-u_t^{*,2}|^2dt\le 0,
\end{align*}
which shows that $u_t^{*,1}=u_t^{*,2}$.

\section*{Acknowledgments}

\qquad The authors acknowledge the financial support from the National Science Foundation of China (grant no. 11871244).

\bibliography{main}

\end{document}